\documentclass[10pt]{article}
\usepackage{amsmath}
\usepackage{amscd}
\usepackage{amsfonts}
\usepackage{amssymb}
\usepackage{amsthm}
\usepackage{mathrsfs}
\usepackage{latexsym}
\usepackage{graphicx}
\usepackage{times}
\usepackage{natbib}
\usepackage{fullpage}
\usepackage{setspace}
\newcommand{\mbs}[1]{\boldsymbol{#1}}
\title{Logistic regression with unknown sizes}
\author{Wei Zhang}
\date{}

\begin{document}
\maketitle
\centerline
{Department of Statistics, University of California, Riverside, CA
, 92521}
\centerline{wxz118@yahoo.com}
\maketitle
\begin{abstract}
Binomial data with unknown sizes often appear in 
biological and medical sciences.
The previous methods either 
use the Poisson approximation
or the quasi-likelihood approach. A full likelihood
approach is proposed by treating unknown sizes
as latent variables.
This approach
simplifies analysis as maximum likelihood estimation can be applied. 
It also facilitates us to gain a lot more insights into 
efficiency loss across models 
and estimation precision within models.  
Simulation assesses the performance of the proposed model. 
An application to the surviving jejunal crypt data is
discussed. The proposed method
is not only competitive with the previous methods, but also 
gives an appropriate
explanation of the inflated variation of expected sizes.
\end{abstract}


\section{Introduction}

Binomial data with unknown sizes often appear in biological and 
medical sciences. For example,
Margolin et al. (1981) studied how the 
number of revertant colonies of Salmonella strain TA98 
changed with the dosage of a chemical agent quinolin. 
Bailer and Piegorsch (2000) reviewed a \textit{C. dubia} 
survival and reproduction toxicity test. 
Trajstman (1989) presented a data set from an experiment of 
 a \textit{M. bovis} subjected to the decontaminants.
Elder (1996) investigated how the times in high heat affect the 
survival of  V79-473 cells.

The jejunal radiation damage is studied extensively in the literature
of medical sciences. Some studies concerned the clinical value of 
some material 
in protecting jejunal crypts
against radiation
(e.g., Goel et al. 2003, Salin et al.  2001 and Khan et al. 1997).
Other studies investigated radiosensitivity of jejunal crypt stem cells
(e.g., Kinashi et al. 1997). 
There are also studies interested in survival of crypt epithelial cells 
in the jejunum of mice exposed to different doses of X-rays
(e.g., Mason et al. 1999).
Table \ref{table:data} presents a surviving  jejunal crypt data
set from an experiment done on $126$ mice 
(Kim and Taylor 1994; Elder et al. 1999).
In such an experiment, each mouse is exposed to a certain dose of
gamma rays, then sacrificed to find out the number of crypts survived.
The total number of crypts before the experiment is unknown, since
the only way to know this number is to sacrifice the mouse
while live mice are required in the experiment.

\begin{table}[ht]
\caption{The jejunal crypt data (dose: the dose of gamma radiation in Gy; 
count: the surviving number of crypts of a mouse).}
\label{table:data}
\begin{center}
\begin{tabular}{ll}
\hline
 dose& count\\
\hline
 6.25 & 76, 96, 73, 81, 81, 87, 77, 75\\
6.50 & 75, 80, 67, 86, 70, 78, 88, 76, 54, 58, 76, 69, 61, 70\\
 6.75 & 66, 51, 48, 48, 57, 45, 59, 49\\
 7.25 & 35, 33, 35, 37, 38, 53, 37, 36, 42, 45, 48, 42, 31, 36, 40, 45, 47,
 38, 40, 35, 27, 35 \\
 7.75 & 19, 18, 25, 19, 19, 18, 21, 18\\
8.00 & 19, 24, 19, 26, 18, 18, 14, 19, 11, 21, 19, 14, 16, 13\\
 8.25 & 19, 19, 19, 16, 12, 16, 12, 13\\
 8.75 & 11, 11, 7, 3, 5, 7, 9, 5, 11, 9, 6, 9,
 7, 5, 10, 7, 11, 9, 7, 11, 5, 12 \\
 9.25 & 6, 3, 5, 6, 4, 6, 5, 3\\
 9.50 & 1, 4, 5, 5, 3, 6, 3, 3, 5, 5, 1, 4, 3, 4\\
\hline
\end{tabular}
\end{center}
\end{table}

Let $y_i$ be the number of surviving crypts in mouse $i$.
It is appropriate to treat $y_i$ as a binomial random variable
with size $n_i$ and surviving probability $p_i$, where
$n_i$ is the total number of crypts in mouse $i$.
The issue of interest is 
to investigate how the surviving probability $p_i$ 
depends on the dose of gamma radiation $x_i$ 
applied to mouse $i$. If the $n_i$ were known, then one could
apply the classical logistic regression (e.g., McCullagh and Nelder 1999).
 Because $n_i$ is unknown,
$y_i$ can also be approximately treated as a Poisson random variable, 
which is a common approach in the literature. 
Such a Poisson approximation is crude when $p_i$ is
moderately large (e.g., Elder et al. 1999). 
By putting additional assumptions on the $n_i$,
  Kim and Taylor (1994) and Elder et al. (1999)
developed a quasi-likelihood approach.
Kim and Taylor (1994) considered that $E(n_i)=m_i$
and $\text{var}\,(n_i)=m_i\nu$ with $m_i$ known and $\nu\geqslant 1$
unknown. Elder et al. (1999) considered estimating $m=E(n_i)$
with $\text{var}\,(n_i)=m(1+\nu m)$ and $\nu \geqslant 0$.

We will assume that each $n_i$ is a Poisson random variable with 
mean $\lambda_i$ and that the $\lambda_i$ arise as a random sample from 
a mixing distribution. In particular,
a gamma distribution will be used in this article. 
By doing this, the requirement of prior knowledge about the $E(n_i)$
in Kim and Taylor (1994) is
 removed. Compared to the quasi-likelihood approach
in Elder et al. (1999), our approach
simplifies analysis as standard techniques, i.e., maximum 
likelihood estimation, can be applied. This model also facilitates us
 to investigate the efficiency loss due to the $n_i$ being unknown and 
being over-dispersed and  how one parameter influence 
  estimation precision of the parameters within a model. Therefore, 
 we can gain much more insights into the problem than previous 
methods.

The proposed method is described in Section 2.
Efficiency losses  are studied in Section 3. 
The estimation precision is investigated in Section 4.
A simulation study is presented in Section 5.
The investigation of the jejunal crypt data is done in Section 6.

\section{The proposed method}
Suppose that the data consist of $r$ pairs of $(y_i,\boldsymbol{x}_i)$, 
$i=1$, 2, $\dots$, $r$, where $\mbs{x}_i$ is the covariate associated
with observation $i$, such that
$y_i\sim \text{Bin}(n_i,p_i)$, $p_i=h(\mbs{x}_i, \mbs{\beta})$.
Note that $h^{-1}$ is a known link function, such as
\[
h(\mbs{x}, \mbs{\beta})=\frac{\exp(\mbs{x}'\mbs{\beta})}
{1+\exp(\mbs{x}'\mbs{\beta})}.
\] 
When each $n_i$ is assumed to be a Poisson random variable with 
mean $\lambda_i$, it is easily shown that
$y_i|\lambda_i\sim \text{Pois}
(\lambda_ih(\mbs{x}_i, \mbs{\beta}))$. 
We will further assume that the $\lambda_i$
arise as a random sample from a gamma density which can be written as
\[
\Gamma^{-1}(\alpha)\eta^{\alpha}\lambda^{\alpha-1}\exp(-\eta\lambda), \lambda \in
(0,\infty),
\]
where $\alpha$ is the shape parameter and $\eta$ is the rate parameter. 
The mean is $\mu=\alpha/\eta$ and the variance is $\sigma^2=\alpha/\eta^2$.
Note that $\alpha=(\mu/\sigma)^2$, that is, $1/\alpha$ is the squared 
coefficient of variation. With
  $(\alpha, \mu)$ used to parameterize gamma densities
 and $\mbs{\theta}=(\mbs{\beta},\mu,\alpha)$,
it is clear that marginally $y_i$ is a negative binomial random variable
 (Anscombe 1949) with density
\[
f(y;\mbs{x},\mbs{\theta})=
\frac{\Gamma(\alpha+y)\alpha^{\alpha}\{\mu h(\mbs{x}, \mbs{\beta})\}^y}
{\Gamma(\alpha)y!\{\alpha+\mu h(\mbs{x}, \mbs{\beta})\}^{\alpha+y}}.
\]
The log likelihood is written as
\begin{equation}
\label{equation:likelihood}
\ell_r(\mbs{\theta})=
\sum_{i=1}^r \log f(y_i;\mbs{x}_i,\mbs{\theta}).
\end{equation}
This looks much like, but is not a special case
of negative binomial regression.
In negative binomial regression, $y|\lambda\sim  \text{Pois}(\lambda)$, 
$\lambda\sim \Gamma(\mu,\mu^2/\alpha)$, and $\mu=\phi(\mbs{x},\mbs{\beta})$ 
for some function $\phi$. 
In the proposed model, 
$y|\lambda\sim \text{Pois}
(\lambda h(\mbs{x}, \mbs{\beta}))$
$\lambda\sim \Gamma(\mu,\mu^2/\alpha)$, and $\mu$ 
is only a parameter and has nothing to do with $\mbs{\beta}$.

Let $\hat{\mbs{\theta}}$ be the maximum likelihood estimator (MLE)
for $\mbs{\theta}$. Asymptotically $\hat{\mbs{\theta}}$
is a multivariate normal random vector with mean $\mbs{\theta}$
and variance-covariance matrix $I_r^{-1}(\mbs{\theta})$
(e.g., Lehmann and Casella 1998, chapter 6), where,
with $\nabla_{\mbs{\beta}}$ being
the gradient with respect to $\mbs{\beta}$,
$I_r(\mbs{\theta})$ is the Fisher information matrix given by
\begin{equation}
\label{equation:itheta}
{I}_r(\mbs{\theta})=
\begin{bmatrix}
\sum_{i=1}^r\frac{\mu
\nabla_{\mbs{\beta}} h(\mbs{x}_{i}, \mbs{\beta})
\nabla_{\mbs{\beta}}'h(\mbs{x}_{i}, \mbs{\beta})
}{h(\mbs{x}_{i}, \mbs{\beta})\{1+\alpha^{-1}\mu h(\mbs{x}_{i}, \mbs{\beta})\}}
&
\sum_{i=1}^r\frac{\nabla_{\mbs{\beta}}h(\mbs{x}_{i}, \mbs{\beta})}
{1+\alpha^{-1}\mu h(\mbs{x}_{i}, \mbs{\beta})}
&0\\
\sum_{i=1}^r\frac{\nabla_{\mbs{\beta}}'
h(\mbs{x}_{i}, \mbs{\beta})}{1+\alpha^{-1}\mu h(\mbs{x}_{i}, \mbs{\beta})}
&\sum_{i=1}^r
\frac{h(\mbs{x}_{i}, \mbs{\beta})}{\mu\{1+\alpha^{-1}\mu h(\mbs{x}_{i}, \mbs{\beta})\}}
&0\\
0&0&
-{E}\left\{\frac{\partial^2 \ell_r(\mbs{\theta})}
{\partial \alpha^2}\right\}
\end{bmatrix}.
\end{equation}

Because of those zero entries in $I_r(\mbs{\theta})$, 
$\alpha$ is orthogonal to the pair of $(\mbs{\beta}$, $\mu)$.
There are several consequences (e.g., Cox and Reid 1987). The asymptotic 
standard errors of $\hat{\mbs{\beta}}$ and $\hat{\mu}$
are not affected by treating $\alpha$ as either known or unknown. The
MLEs of $(\mbs{\beta}$, $\mu)$ and $\alpha$ are asymptotically independent.
The MLEs of $\mbs{\beta}$ and $\mu$ given $\alpha$  
vary only slowly with $\alpha$.

\section{Efficiency loss}
The parameter of interest is $\mbs{\beta}$. 
Within the proposed model, the MLE  $\hat{\mbs{\beta}}$
is asymptotically fully efficient.
If we knew the $n_i$, then a more precise estimation
of $\mbs{\beta}$
is feasible. Two kinds of efficiency losses 
are of interest: 
that originated from the $n_i$ being unknown, and
that from the over-dispersion among the $n_i$. 

If we knew the $n_i$, then the Fisher information matrix of $\mbs{\beta}$  is
\[
{I}_r(\mbs{\beta}|\{n_i\}_{i=1}^r)=
\sum_{i=1}^r\frac{n_i\nabla_{\mbs{\beta}}h(\mbs{x}_i, \mbs{\beta})
\nabla_{\mbs{\beta}}'h(\mbs{x}_i, \mbs{\beta})}
{h(\mbs{x}_i, \mbs{\beta})\{1-h(\mbs{x}_i, \mbs{\beta})\}}.
\]
When the $n_i$ are unknown and arise
from  a Poisson  distribution  with mean $\mu$,
the Fisher information matrix is 
\begin{equation}
\label{equation:ibetamu}
{I}_r(\mbs{\beta},\mu)=
\begin{bmatrix}
\sum_{i=1}^r\frac{\mu
\nabla_{\mbs{\beta}} h(\mbs{x}_{i}, \mbs{\beta})
\nabla_{\mbs{\beta}}'h(\mbs{x}_{i}, \mbs{\beta})
}{h(\mbs{x}_{i}, \mbs{\beta})}
&
\sum_{i=1}^r\nabla_{\mbs{\beta}}h(\mbs{x}_{i}, \mbs{\beta})\\
\sum_{i=1}^r\nabla_{\mbs{\beta}}'
h(\mbs{x}_{i}, \mbs{\beta})
&\sum_{i=1}^r
\frac{h(\mbs{x}_{i}, \mbs{\beta})}{\mu}
\end{bmatrix},
\end{equation}
which is obtained from ${I}_r(\mbs{\theta})$ in 
(\ref{equation:itheta}) by letting $\alpha=\infty$.

A sensible approach to evaluate the efficiency loss
is to average ${I}_r(\mbs{\beta}|\{n_i\}_{i=1}^r)$
by treating  the $n_i$ as a sample from  a Poisson  distribution  
with mean $\mu$,
i.e., to consider
\begin{equation}
\label{equation:ibetagivenmu}
{I}_r(\mbs{\beta}|\mu)=
E\{{I}_r(\mbs{\beta}|\{n_i\}_{i=1}^r)\}=
\sum_{i=1}^r\frac{\mu\nabla_{\mbs{\beta}}h(\mbs{x}_i, \mbs{\beta})
\nabla_{\mbs{\beta}}'h(\mbs{x}_i, \mbs{\beta})}
{h(\mbs{x}_i, \mbs{\beta})\{1-h(\mbs{x}_i, \mbs{\beta})\}}.
\end{equation}

A numeric experiment is used to investigate the efficiency losses,
in which there is a single covariate $x$ and the parameter
of interest is the slope $\beta_1$ while the intercept $\beta_0$
is fixed to be one. A $2^4$ design is considered,  i.e., 
\begin{equation}
\underbrace{\{\mathcal{X}_1,\mathcal{X}_2\}}_{x}
\times \underbrace{\{1,2\}}_{\beta_1}
\times \underbrace{\{100,300\}}_{\mu}
\times \underbrace{\{25,49\}}_{\alpha},
\end{equation}
where $\mathcal{X}_1$
is the set of integers in $[-5,5]$, and $\mathcal{X}_2$ is
$\{-0.63$, $1.59$, $-3.01$,$ -6.85$, $-4.97$, $ 1.86$, 
$-7.54$, $-3.45$, $-4.45$, $-1.87$, $6.49\}$,
a set of $11$ normal random variables with mean $0$ and variance $25$.
It is assumed that  the number of replications
is identical over each $x$ value. 

Let $\rho$ be the ratio of the asymptotic standard deviation
of $\beta_1$ calculated from ${I}_r(\mbs{\beta}|\mu)$ in
(\ref{equation:ibetagivenmu})
and that from ${I}_r(\mbs{\beta},\mu)$ in (\ref{equation:ibetamu}).
Let $\gamma$ be the ratio of the asymptotic standard deviation
of $\beta_1$ calculated from ${I}_r(\mbs{\beta},\mu)$ in 
(\ref{equation:ibetamu})
and that from ${I}_r(\mbs{\theta})$ in (\ref{equation:itheta}).
This means that $\rho\gamma$ is the ratio of the asymptotic standard deviation
of $\beta_1$ calculated from ${I}_r(\mbs{\beta}|\mu)$ in 
(\ref{equation:ibetagivenmu})
and that from ${I}_r(\mbs{\theta})$ in (\ref{equation:itheta}).
Table \ref{table:effloss} presents these efficiency loss measures.
 The ranges of
 $\rho$, $\gamma$ and $\rho \gamma$  are given by 
$(0.706,0.786)$, $(0.732, 0.941)$ and
 $(0.517,0.740)$, respectively. 

\begin{figure}
\centering
\includegraphics[angle=270,width=12cm]{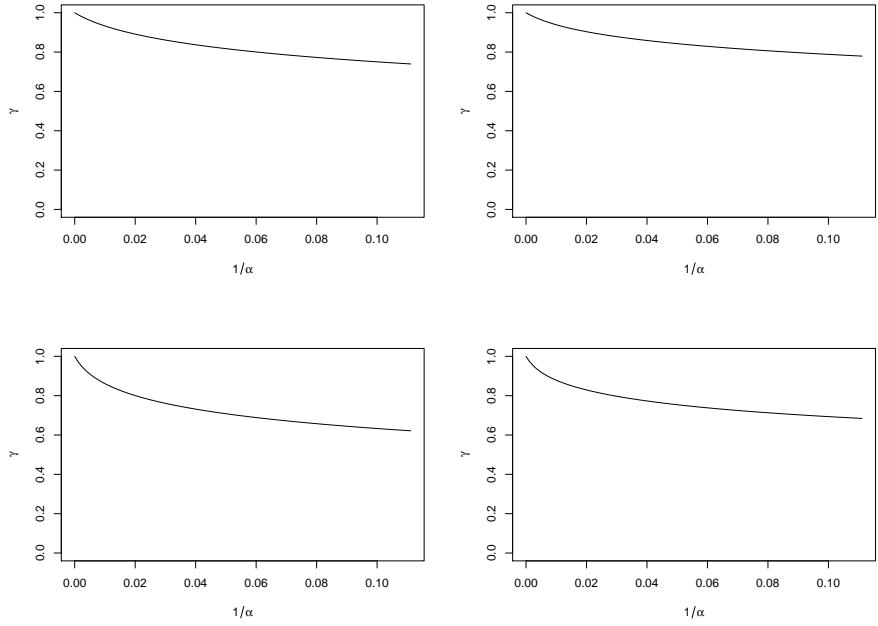}
\caption{The efficiency loss measure $\gamma$ varies 
with $\alpha$: $\beta_0=1$, $x\in\mathcal{X}_1$,
$(\beta_1,\mu)=(1,100)$ (upper left panel),
$(2,100)$ (upper right), $(1,300)$ (lower left), 
and $(2,300)$ (lower right). 
 }
\label{figure:efficiency}
\end{figure}

Figure \ref{figure:efficiency} shows how 
the efficiency loss measure $\gamma$ changes 
when  $\alpha$ varies continuously.
As $\alpha$ increases, the gamma distribution tends to be degenerated, 
and  the efficiency loss from the over-dispersion among the $n_i$  decreases.
In the four panels, the efficiency loss is at most 0.621,
which indicates that the efficiency loss from the over-dispersion among 
the $n_i$ is  small.

\begin{table}[ht]
\caption{The efficiency loss measures over $16$ settings.}
\label{table:effloss}
\begin{center}
\begin{tabular}{rrrrrrr}
\hline
setting & $\beta_1$ & $\mu$ & $\alpha$ & $\rho$ & $\gamma$ & $\rho \gamma$ \\
\hline
1 & 1 & 100 & 25 & 0.706 & 0.837 & 0.591 \\
2 & 2 & 100 & 25 & 0.747 & 0.859 & 0.642 \\
3 & 1 & 300 & 25 & 0.706 & 0.732 & 0.517 \\
4 & 2 & 300 & 25 & 0.747 & 0.773 & 0.578 \\
5 & 1 & 100 & 49 & 0.706 & 0.890 & 0.629 \\
6 & 2 & 100 & 49 & 0.747 & 0.902 & 0.674 \\
7 & 1 & 300 & 49 & 0.706 & 0.799 & 0.564 \\
8 & 2 & 300 & 49 & 0.747 & 0.828 & 0.618 \\
9 & 1 & 100 & 25 & 0.729 & 0.858 & 0.625 \\
10 & 2 & 100 & 25 & 0.786 & 0.902 & 0.709 \\
11 & 1 & 300 & 25 & 0.729 & 0.765 & 0.558 \\
12 & 2 & 300 & 25 & 0.786 & 0.816 & 0.641 \\
13 & 1 & 100 & 49 & 0.729 & 0.904 & 0.659 \\
14 & 2 & 100 & 49 & 0.786 & 0.941 & 0.740 \\
15 & 1 & 300 & 49 & 0.729 & 0.824 & 0.601 \\
16 & 2 & 300 & 49 & 0.786 & 0.871 & 0.685 \\
\hline
\end{tabular}
\end{center}
\end{table}

\section{Estimation precision of $\mbs{\beta}$ and $\mu$}
 
For a Poisson sample of size $r$ with mean $\mu$, the information of $\mu$
is given by ${I}_r(\mu)=r/\mu$. When $\mu$ increases, the estimation of $\mu$ 
becomes less precise.
On the contrary, a large $\mu$ serves a good purpose for the estimation
of  $\mbs{\beta}$ in a logistic regression model, which can be
clearly seen from ${I}_r(\mbs{\beta}|\mu)$ in (\ref{equation:ibetagivenmu}).

When the $n_i$ are unknown and arise from the a Poisson distribution with mean
$\mu$, although $\mbs{\beta}$ and $\mu$ are not orthogonal, we will show that
a large $\mu$ will lead to a more precise estimation of $\mbs{\beta}$
 but a less precise estimation of $\mu$.
To this end, we will partition the
asymptotic variance-covariance matrix $V$, i.e., 
the inverse matrix of ${I}_r(\mbs{\beta},\mu)$ in (\ref{equation:ibetamu}),
 into a $2\times 2$ block form, where, 
$V=(V_{ij})$, and
\begin{align*}
V_{11}&=\frac{1}{\mu}\Bigl\{\sum_{i=1}^r\frac{
\nabla_{\mbs{\beta}} h(\mbs{x}_{i}, \mbs{\beta})
\nabla_{\mbs{\beta}}'h(\mbs{x}_{i}, \mbs{\beta})
}{h(\mbs{x}_{i}, \mbs{\beta})}-
\frac{\sum_{i=1}^r
\nabla_{\mbs{\beta}} h(\mbs{x}_{i}, \mbs{\beta})
\sum_{i=1}^r\nabla_{\mbs{\beta}}' h(\mbs{x}_{i}, \mbs{\beta})}
{\sum_{i=1}^rh(\mbs{x}_{i}, \mbs{\beta})}\Bigr\}^{-1},\\
V_{22}&=\mu\Bigl[\sum_{i=1}^rh(\mbs{x}_{i}, \mbs{\beta})-
\sum_{i=1}^r\nabla_{\mbs{\beta}}'h(\mbs{x}_{i}, \mbs{\beta})
\Bigl(\sum_{i=1}^r\frac{
\nabla_{\mbs{\beta}} h(\mbs{x}_{i}, \mbs{\beta})
\nabla_{\mbs{\beta}}'h(\mbs{x}_{i}, \mbs{\beta})
}{h(\mbs{x}_{i}, \mbs{\beta})}\Bigr)^{-1}\sum_{i=1}^r
\nabla_{\mbs{\beta}} h(\mbs{x}_{i}, \mbs{\beta})\Bigr]^{-1}.
\end{align*}
The diagonal entries of $\mu V_{11}$ and $V_{22}/\mu$
are nonnegative. The variance of each component in $\mbs{\beta}$
 is a nonincreasing 
function of $\mu$,
while that of $\mu$ is a nondecreasing function of $\mu$.

\begin{figure}
\centering
\includegraphics[angle=270,width=12cm]{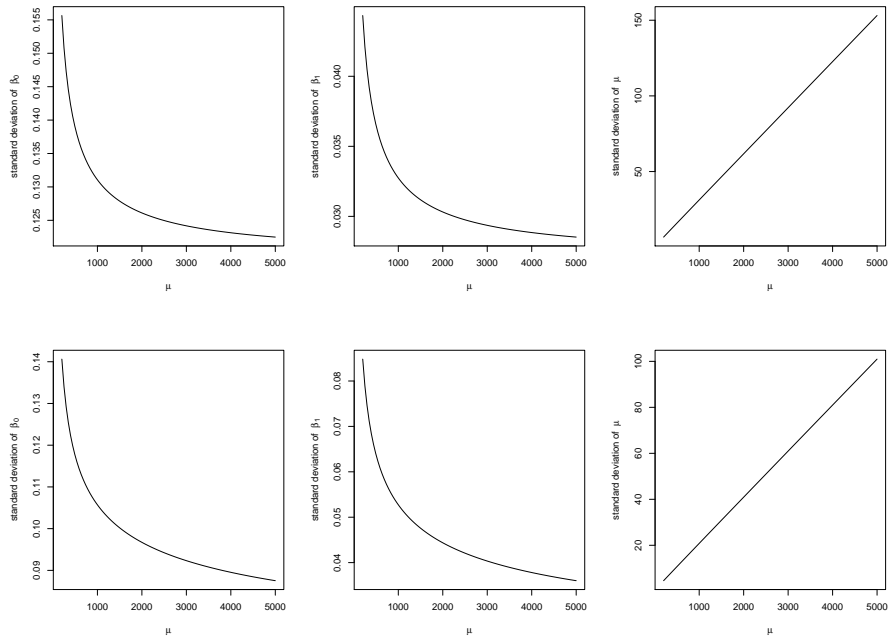}
\caption{The standard deviations of $\mbs{\beta}$ and $\mu$
change with respect to $\mu$: $\beta_0=1$, $x\in\mathcal{X}_1$, 
$10$ replications at each $x$,  
$(\beta_1, \alpha)=(1, 25)$ (top panels)
and $(2, 49)$ (bottom panels).}
\label{figure:varmu}
\end{figure}

When the $n_i$ are over-dispersed, we also conjecture that 
a large $\mu$ will have the same
 effects on the estimation of 
$\mu$ and $\mbs{\beta}$ as those in the simple Poisson model.
Figure \ref{figure:varmu}  is
a numerical illustration. 

\section{Simulation}

The simulation with the same $16$ settings
as the efficiency loss study and $1000$ samples is reported in
Table \ref{table:simres}.  
The number of replications is $10$ over each $x$ value. The 
bias and mean square error are pretty small. All coverage probabilities
of $95\%$ confidence interval achieve their nominal value 0.95.

\begin{table}[ht]
\caption{Simulation results (the nominal value of coverage
probability is $0.95$).}
\label{table:simres}
\begin{center}
\begin{tabular}{rrrrrcc}
\hline
setting& $\beta_1$ & $\mu$ & $\alpha$ & bias & mean square error& 
coverage probability \\
\hline
1 & 1 & 100 & 25 & 0.005 & 0.003 & 0.956 \\
2 & 2 & 100 & 25 & 0.007 & 0.015 & 0.956 \\
3 & 1 & 300 & 25 & 0.003 & 0.002 & 0.947 \\
4 & 2 & 300 & 25 & 0.008 & 0.007 & 0.948 \\
5 & 1 & 100 & 49 & 0.005 & 0.002 & 0.951 \\
6 & 2 & 100 & 49 & 0.004 & 0.013 & 0.927 \\
7 & 1 & 300 & 49 & 0.002 & 0.001 & 0.951 \\
8 & 2 & 300 & 49 & 0.004 & 0.005 & 0.945 \\
9 & 1 & 100 & 25 & 0.000 & 0.003 & 0.949 \\
10 & 2 & 100 & 25 & 0.014 & 0.015 & 0.948 \\
11 & 1 & 300 & 25 & 0.001 & 0.001 & 0.946 \\
12 & 2 & 300 & 25 & $-$0.001 & 0.007 & 0.952 \\
13 & 1 & 100 & 49 & 0.001 & 0.002 & 0.948 \\
14 & 2 & 100 & 49 & 0.014 & 0.013 & 0.949 \\
15 & 1 & 300 & 49 & 0.001 & 0.001 & 0.947 \\
16 & 2 & 300 & 49 & 0.004 & 0.006 & 0.945 \\
\hline
\end{tabular}
\end{center}
\end{table}

\section{Example}

For the jejunal crypt data in Table \ref{table:data}, 
it is assumed that the surviving probabilities $p_i$ satisfy 
$\log\{p_i/(1-p_i)\}=\beta_0+\beta_1x_i$ for all $i$ (e.g., 
Kim and Taylor 1994 and Elder et al. 1999).
The R function \texttt{optim} is used to maximize the likelihood function 
in (\ref{equation:likelihood}).
The estimates of $\mbs{\beta}$ and $\mu$ are stable, but that
of $\alpha$ varies a lot. The estimate $\hat{\alpha}$ is 
$3121.834$ when its initial value is $20$, but 
becomes $6070.602$ when its initial value is $200$. 
The Hessian matrix is found to be nearly singular,
which implies that the variance of $\alpha$ is huge.

By the likelihood ratio test, 
we would like to assume that the $n_i$ arise from a single Poisson
distribution.
The results are reported in Table \ref {table:eg}, 
which also shows the estimates and their standard errors
using the logistic regression (with $n_i=160$), 
the quasi-likelihood approaches in 
Kim and Taylor (1994) (with $E(n_i)=160$) and Elder et al. (1999).
Our estimates and standard errors are very close to those in Elder 
et al. (1999).
All the estimates of previous methods fall into our $95\%$ 
confidence intervals: $(5.207,\, 8.203)$ for $\beta_0$,
$(-1.248,\, -1.000)$ for $\beta_1$ and 
$(103.4,\, 289.0)$ for $\mu$. The standard errors of $\mbs{\beta}$
are pretty small, while that of $\mu$ is quite large. 
Elder et al. (1999) conjectured that the variance inflation of $\mu$ 
 is due to the data structure, i.e., 
there is no zero dose. By the estimation precision study in Section 4, 
we can not only give a more accurate explanation of the large standard error 
of $\mu$, but also explains the small standard errors of $\mbs{\beta}$. 
 Since $\hat{\mu}$ is as large as $196.2$,
the standard error of $\mu$ is large, while those of
$\mbs{\beta}$ are small.

\begin{table}[ht]
\caption{The jejunal crypt data results from the proposed and
previous approaches 
(logistic regression and Kim's method fix $n_i$ and $E(n_i)$ at 160, 
 respectively; Kim's and Elder's quasi-likelihood method of moments estimates
come from Elder et al. (1999)).}
\label{table:eg}
\begin{center}
\begin{tabular}{c|ccccc}
\hline
& \multicolumn{4}{c}{estimate (standard error)}\\
\cline{2-5}
&logistic &Kim's&
 Elder's&proposed\\
\hline
$\beta_0$ & 7.432 (0.175)&7.410 (0.191) &6.727 (0.725)&6.705 (0.764)
\\
$\beta_1$ & $-$1.185 (0.024) &$-$1.183 (0.026)& $-$1.126 (0.061)&$-$
1.124 (0.063)\\
$\mu$ &--- &--- & 194.7 (43.4)&196.2 (47.4) \\
\hline
\end{tabular}
\end{center}
\end{table}

\section{Discussion}

One may consider estimating the $n_i$ and then apply
the logistic regression. 
There is a lot of literature about estimating binomial size $n$ 
under the condition that
$p$ is either known or unknown. If $p$ is unknown,
then it is usually treated as a nuisance parameter 
(e.g., Draper and Guttman, 1971; Caroll and Lombard, 1985).  
Unlike the studies in the literature,
the $p_i$ depend on covariates, and many
$n_i$ need to be estimated. 
Such a two-stage approach also makes analysis unnecessarily
 more complicated. The proposed approach estimates all parameters
in a seamless fashion by treating the means of the $n_i$ as nuisance
parameters and integrating them out.

\section*{Appendix: The Fisher information matrix}
Let $\ell(\mbs{\theta})=\log f(y;\mbs{x},\mbs{\theta})$
and $\varPsi(\alpha)=\log\Gamma(\alpha)$.
The first order derivatives are
\begin{align*}
\frac{\partial \ell}{\partial\alpha}&=
\varPsi'(\alpha+y)-\varPsi'(\alpha)+\log\alpha+1-\log
\{\alpha+\mu h(\mbs{x}, \mbs{\beta})\}-
\frac{\alpha+y}{\alpha+\mu h(\mbs{x}, \mbs{\beta})},
\\
\frac{\partial \ell}{\partial\mu}&=
\frac{y}{\mu}-\frac{\alpha+y}
{\alpha+\mu h(\mbs{x}, \mbs{\beta})}
h(\mbs{x}, \mbs{\beta}),
\\
\frac{\partial \ell}{\partial\mbs{\beta}}&=
\frac{y}{
 h(\mbs{x}, \mbs{\beta})}\nabla_{\mbs{\beta}} h(\mbs{x}, \mbs{\beta})
-\frac{(\alpha+y)\mu}{\alpha+\mu h(\mbs{x}, \mbs{\beta})}
 \nabla_{\mbs{\beta}} h(\mbs{x}, \mbs{\beta}).
\end{align*}
The second order derivatives are
\begin{align*}
\frac{\partial^2 \ell}{\partial \alpha^2}&
=\varPsi''(\alpha+y)-\varPsi''(\alpha)+\frac{1}{\alpha}
-\frac{1}{\alpha+\mu h(\mbs{x}, \mbs{\beta})}
+\frac{y-\mu h(\mbs{x}, \mbs{\beta})}
{\{\alpha+\mu h(\mbs{x}, \mbs{\beta})\}^2},\\
\frac{\partial^2 \ell}{\partial \mu^2}&
=-\frac{y}{\mu^2}+\frac{(\alpha+y) h^2(\mbs{x}, \mbs{\beta})}
{\{\alpha+\mu h(\mbs{x}, \mbs{\beta})\}^2},\\
\frac{\partial^2 \ell}{\partial \mbs{\beta}\partial\mbs{\beta}'}&
=\left\{
\frac{y}{h(\mbs{x}, \mbs{\beta})}-\frac{(\alpha+y)\mu}{\alpha+\mu h(\mbs{x}, \mbs{\beta})}
\right\}\frac{\partial^2 h(\mbs{x}, \mbs{\beta})}
{\partial \mbs{\beta}\partial\mbs{\beta}'}\\
&\phantom{=}+
\left\{
\frac{(\alpha+y)\mu^2}{\{\alpha+\mu h(\mbs{x}, \mbs{\beta})\}^2}
-\frac{y}{h^2(\mbs{x}, \mbs{\beta})}
\right\}\nabla_{\mbs{\beta}}h(\mbs{x}, \mbs{\beta})
\nabla_{\mbs{\beta}}'h(\mbs{x}, \mbs{\beta}),\\
\frac{\partial^2 l}{\partial \alpha\partial\mu}&
=\frac{h(\mbs{x}, \mbs{\beta})\{y-\mu h(\mbs{x}, \mbs{\beta})\}}
{\{\alpha+\mu h(\mbs{x}, \mbs{\beta})\}^2},
\\
\frac{\partial^2 \ell}{\partial \alpha\partial\mbs{\beta}}&
=\frac{\mu\{y-\mu h(\mbs{x}, \mbs{\beta})\}}
{\{\alpha+\mu h(\mbs{x}, \mbs{\beta})\}^2}
\nabla_{\mbs{\beta}} h(\mbs{x}, \mbs{\beta}),
\\
\frac{\partial^2 \ell}{\partial \mu\partial\mbs{\beta}}
&=\frac{-(\alpha+y)\alpha}{\{\alpha+\mu h(\mbs{x}, \mbs{\beta})\}^2}
\nabla_{\mbs{\beta}} h(\mbs{x}, \mbs{\beta}).
\end{align*}
By taking negative expectation with respect to $f(y;\mbs{x},\mbs{\theta})$,
one obtains the Fisher information matrix. Note that
$E(y)=\mu h(\mbs{x}, \mbs{\beta} )$.

\nocite{*}
\bibliography{ref}
\bibliographystyle{apalike}

\end{document}